\documentclass[reqno,tbtags,11pt]{amsproc}
\usepackage{amssymb,ifthen}

\def\YEAR{01}
\def\VOLUME{35}
\def\ISSUE{3}
\def\MONTHES{July--September}

\newcounter{pppage}
\newsavebox{\PAGES}
\newsavebox{\EPAGES}
\newsavebox{\UDC}
\newsavebox{\PLENUM}
\newsavebox{\DATA}
\newsavebox{\TRANSL}
\newsavebox{\EDITOR}

\let\rom=\textup
\def\rmm#1{\kern-\mathsurround$#1$\kern-\mathsurround}
\def\cit#1{\rmm{[#1]}}
\def\citeoff{\let\cite=\cit}
\def\eq#1{\rmm{(#1)}}

\def\CC{\mathbb{C}}

\def\ZZ{\mathbb{Z}}

\catcode`\@=11

\def\@biblabel#1{#1.}
\def\@cite#1#2{{%
 \m@th\upshape\mdseries[{\mdseries #1}{\if@tempswa, #2\fi}]}}
\@ifundefined{cite }{%
  \expandafter\let\csname cite \endcsname\cite
  \edef\cite{\@nx\protect\@xp\@nx\csname cite \endcsname}%
}{}
\def\footnoterule{\kern11pt\kern-\footnotesep}
\def\th@plain{%
  \let\thm@indent\indent
  \thm@headfont{\bfseries}%
  \let\thmhead\thmhead@plain \let\swappedhead\swappedhead@plain
  \thm@preskip.25\baselineskip\@plus.1\baselineskip
                                    \@minus.1\baselineskip
  \thm@postskip\thm@preskip
  \itshape
}
\def\th@definition{%
  \let\thm@indent\indent
  \thm@headfont{\bfseries}
  \let\thmhead\thmhead@plain \let\swappedhead\swappedhead@plain
  \thm@preskip.25\baselineskip\@plus.1\baselineskip
                                    \@minus.1\baselineskip
  \thm@postskip\thm@preskip
  \upshape
}
\def\th@remark{%
  \let\thm@indent\indent
  \thm@headfont{\bfseries}
  \let\thmhead\thmhead@plain \let\swappedhead\swappedhead@plain
  \thm@preskip.25\baselineskip\@plus.1\baselineskip
                                    \@minus.1\baselineskip
  \thm@postskip\thm@preskip
  \upshape
}
\catcode`\@=13

\def\pages#1{\savebox{\PAGES}{\footnotesize #1}}

\def\translator{\bigskip\hbox to \textwidth{\hfil Translated by \usebox{\TRANSL}}}

\def\title#1#2#3#4#5{
\centerline{\Large\bf #1}
\centerline{\Large\bf #2
\ifthenelse{\equal{#3}{}}{}{\rm*}}
\def\thefootnote{}\footnotesep=10pt
\ifthenelse{\equal{#3}{}}{\protect\footnotetext{\capt{-\parindent}%
{-12.3pt}{\rule[-3pt]{\textwidth}{0.4pt}}#4. Published in
Funktsio\-nal\smash{\mz}nyi
Ana\-liz i Ego Pri\-lo\-zhe\-niya, Vol.~\VOLUME, No.~\ISSUE,
pp.~\usebox{\PAGES},
\MONTHES, 20\YEAR. Original article submitted #5.}}%
{\protect\footnotetext{\capt{-\parindent}{-12.3pt}{\rule[-3pt]{10pc}{0.4pt}}%
\ifthenelse{\equal{#3}{}}{}{*}#3.\rule[-12pt]{0pt}{11pt}\newline
\capt{0pt}{-12.3pt}{\rule[-3pt]{\textwidth}{0.4pt}}\kern\parindent#4.
Published in Funktsio\-nal\smash{\mz}nyi
Ana\-liz i Ego Pri\-lo\-zhe\-niya, Vol.~\VOLUME, No.~\ISSUE,
pp.~\usebox{\PAGES},
\MONTHES, 20\YEAR. Original article submitted #5.}}\vskip4pt}

\def\author#1{\centerline{#1}}

\def\subhead#1{\vskip3pt\textbf{\bf #1.}}

\def\le{\leqslant}
\def\ge{\geqslant}

\def\={\,{=}\,}

\def\a{\alpha}

\def\ga{\gamma}

\def\la{\lambda}
\def\e{\varepsilon}

\def\ov{\overline}

\def\â.¥.{â.{\kern1.6pt}¥.}
\def\¯.¢.{¯.{\kern1.6pt}¢.}

\def\C{\mathbb{C}\kern0.5pt}
\def\R{\mathbb{R}\kern0.5pt}

\DeclareMathOperator\rank{rank}

\newcommand{\op}[1]{\operatorname{#1}}

\font\tensscr=rsfs10.pk scaled 1000
\font\sevensscr=rsfs7.pk scaled 1000
\font\fivesscr=rsfs5.pk scaled 1000
\def\mathCal#1{\mathchoice{\hbox{\tensscr #1}}{\hbox{\tensscr #1}}%
              {\hbox{\sevensscr #1}}{\hbox{\fivesscr #1}}}

\newenvironment{tm}[1]{\ifdim\lastskip<2.1pt\vskip%
                       -\lastskip\else\fi\vskip2pt\textbf{#1.}\em}{\vskip2pt}
\newenvironment{pf}[1]{\ifdim\lastskip<2.1pt\vskip%
                       -\lastskip\else\fi\vskip2pt\textbf{#1.}}{\vskip2pt}
\newenvironment{ab}{\leftskip20pt\rightskip20pt\baselineskip10.5pt\small}{\vskip2pt}

\def\cA{{\mathCal{A}}}

\def\cM{{\mathCal{M}}}

\def\cL{{\mathCal{L}}}

\def\ovl#1#2#3{\rlap{$#3$}\kern#1pt\ov{\phantom{#3}%
                                 \kern-#1pt\kern-#2pt}\kern#2pt}
\def\mz{$'$}
\def\wt{\widetilde}

\def\itemi{\par\parindent19.8pt{\rm(i)}\kern4pt}
\def\itemii{\par\parindent17pt{\rm(ii)}\kern4pt}
\def\itemiii{\par\parindent14pt{\rm(iii)}\kern4pt}
\def\itemiv{\par\parindent14.5pt{\rm(iv)}\kern4pt}
\def\itemv{\par\parindent17pt{\rm(v)}\kern4pt}
\def\itemvi{\par\parindent8.5pt{\rm(vi)}\kern4pt}

\def\scirc{\kern2pt\lower-2pt\hbox{$_\circ$}\kern2pt}
\let\rom\textup

\def\capt#1#2#3{\rlap{\kern #1\smash{\lower #2\hbox{#3}}}}
\def\pict#1#2#3#4#5{\rlap{\kern #1\smash{\lower #2\hbox{\includegraphics[bb= 0 0 #3 #4]{#5}}}}}
\def\pictPS#1#2#3{\rlap{\kern #1\smash{\lower #2\hbox{\includegraphics{#3}}}}}

\def\afline#1#2{}{}
\def\ft#1{\def\footnoterule{\kern-.4pt
        \hrule width 10pc\kern11pt\kern-\footnotesep}%
$^*$\footnote{$^*$#1}}
\def\ftt#1{\def\footnoterule{\kern-.4pt
        \hrule width 10pc\kern11pt\kern-\footnotesep}%
$^{**}$\footnote{$^{**}$#1}}

\begin{document}

\pages{60--72}
\title{The Fermion Model of Representations}
{of Affine Krichever--Novikov Algebras}
{Supported by RFBR grant No. 99-01-00198}
{V.~A.~Steklov Mathematical Institute and Independent University
of Moscow}
{May 11, 2000}

\author{{\bf O.~K.~Sheinman}}

\numberwithin{equation}{section}

\begin{ab}
To a generic holomorphic vector bundle on an algebraic curve
and an irreducible finite-dimen\-sional representation of a
semisimple Lie algebra, we assign a representation of the
corresponding affine Krichever--Novikov algebra in the space
of semi-infinite exterior forms. It is shown that equivalent
pairs of data give rise to equivalent representations and
vice versa.
\end{ab}

\section{Introduction}

Affine Krichever--Novikov algebras appeared in
\cite{rKNFa,rKNFb} as a generalization of affine Kac--Moody
algebras that is related to a compact Riemann surface with
two distinguished points. They belong to the new class of
\textit{quasigraded algebras} (see the definition in \S2),
whose structure and representations are at present at the
initial stage of investigation. At the same time, the notion
of a quasigraded algebra is of importance and enables one to
generalize many properties of the well-studied class of
graded algebras, in particular, the most important property
that one can construct representations generated by the
vacuum vector. The first constructions of representations of
affine Krichever--Novikov algebras (Verma modules and
irreducible modules) were suggested in
\cite{rSheb,rShhw,rShns}, and, in the case of several
distinguished points, in~\cite{rSSkz}. These constructions
are of abstract algebraic nature. The problem of finding
geometric objects on which these algebras act in a natural
way remained open. In this paper we present such an object,
namely, function spaces related in a special way to
holomorphic vector bundles over the corresponding Riemann
surface (these spaces were introduced earlier in
\cite{rKNFd,rKNU} in connection with the solution of soliton
equations).

Starting from these spaces, we construct representations
generated by highest vectors by\break means of the known
construction of an infinite fermion representation (the
wedge representation)\break
\cite{rKaPet,rKNFa,rKNFb,rSLc,rFeFu}; see also
\cite{rKaRa,rKac}. To this end, we use some bases in the
function spaces under consideration; these bases were
introduced in \cite{rKNRm} and generalize the well-known
Krichever--Novikov bases in function and tensor spaces on a
Riemann surface with two distinguished
points~\cite{rKNFa,rKNFb,rKNFc}. In \S2 we give a
construction of these bases and slightly extend the
exposition in \cite{rKNRm}.

The main theorem of the paper, Theorem~3.2, establishes
conditions for the equivalence of fermion representations,
and the most substantial condition is the equivalence of
bundles generating these representations.

The appearance of bundles in representation theory of
Krichever--Novikov algebras is by no means occasional. It
was shown in \cite{rShea,rSha} that the orbit space of the
coadjoint representation of such an algebra coincides with
the space of equivalence classes of finite-dimensional
irreducible representations of the fundamental group of a
punctured Riemann surface. In fact, these are just the data
that arose earlier in the classification of
Narasimhan--Seshadri modules of stable bundles on a
(compact) Riemann surface. These data are also closely
related \cite{rShSD} to the Hitchin construction of Higgs
bundles.

Krichever and Novikov~\cite{rKNRm} arrived at the idea of
constructing representations by using algebraic curves for
other reasons. The quasigraded structure on the
Krichever--Novikov algebras and modules over these algebras
enables one to treat the representation theory of these
algebras as a part of the theory of difference operators,
where the theory of commuting difference operators developed
in the papers of these authors \cite{rKr,rKNRs,rKNRm} plays
a substantial role and readily leads to holomorphic bundles
over curves. In fact, they had in mind the relationship with
the theory of solitons already when writing the papers
\cite{rKNFa,rKNFb,rKNFc}. Note that the deep connection
between the representation theory of affine
Krichever--Novikov algebras and the theory of commutative
rings of difference operators is revealed, in particular, in
the proof of Theorem~3.1 in the present paper.

This paper appeared as the result of my discussions with
I.~M.~Krichever of all specific features of the construction
in question, and the author is deeply indebted to
I.~M.~Krichever for these discussions. Moreover,
I.~M.~Krichever told the author about the results of the
paper \cite{rKNRm} long before the publication, expressed a
lot of ideas about the application of these results in the
posed problem of representation theory, and patiently
discussed versions of their realization. Of these
suggestions, the most significant was to pass from the space
of sections of a vector bundle (on which, as is known, there
is no action of matrix algebras in general) to the
corresponding space of vector functions (see \S2), where the
desired action exists.

\section{Krichever--Novikov Bases in Spaces of Sections
of Holomorphic Bundles}

\subhead{(a) Krichever--Novikov algebras} Let $\Sigma$ be a
compact algebraic curve over $\CC$ with two distinguished
points $P_{\pm}$, let $\cA\,(\Sigma,P_\pm)$ be the algebra
of meromorphic functions on $\Sigma$ regular outside
$P_\pm$, and let $\mathfrak{g}$ be a complex semisimple Lie
algebra. Then
\begin{equation}
\hat{\mathfrak{g}}
=\mathfrak{g}\otimes_\CC\cA\,(\Sigma,P_\pm)\oplus
\CC c\label{alg}
\end{equation}
is called a \textit{Krichever--Novikov algebra of affine
type} \cite{rKNFa,rShea,rSha} with one-dimensional center
generated by the element $c$. Elements of
$\hat{\mathfrak{g}}$ will be denoted by $\wt{X}=X+ac$, where
$X\in\mathfrak{g}\otimes_\CC\cA\,(\Sigma,P_\pm)$ and
$a\in\CC$. The bracket on $\hat{\mathfrak{g}}$ is defined by
$$
[x\otimes A,y\otimes B]=[x,y]\otimes AB
+\ga(x\otimes A,y\otimes B)c,\qquad [x\otimes A,c]=0,
$$
where $\ga$ is the cocycle defined by the relation
\begin{equation}
\ga(x\otimes A,y\otimes B)
=(x,y)\op{res}_{P_+}(A\,dB).\label{cocycle}
\end{equation}

As was mentioned in the introduction, the algebra
$\hat{\mathfrak{g}}$ carries a remarkable quasigraded
structure.

\begin{pf}{Definition 2.1}
(a) Let $\cL$ be a Lie algebra or an associative algebra
that admits a decomposition $\cL=\bigoplus_{n\in\ZZ}\cL_n$
into a direct sum of finite-dimensional subspaces $\cL_n$.
The algebra $\cL$ is said to be \textit{quasigraded}
(\textit{almost graded, graded in the extended sense}) if
$\dim\cL_n<\infty$ and there exist constants $R$ and $S$
such that
\begin{equation}
\cL_n\cdot\cL_m\subseteq
\bigoplus_{h=n+m-R}^{n+m+S}
\cL_h\qquad\forall n,m\in\ZZ.\label{eaga}
\end{equation}
The elements of the subspaces $\cL_n$ are called
\textit{homogeneous elements of degree~$n$}.

(b) Let $\cL$ be a quasigraded Lie algebra or an associative
algebra, and let $\cM$ be an $\cL$-module that admits a
decomposition $\cM=\bigoplus_{n\in\ZZ}\cM_n$ into a direct
sum of subspaces. The module $\cM\,$ is said to be
\textit{quasigraded} (\textit{almost graded, graded in the
extended sense}) if $\dim\cM_n<\infty$ and there exist
constants $R'$ and $S'$ such that
\begin{equation}
\cL_m\cdot\cM_n\subseteq\bigoplus_{h=n+m-R'}^{n+m+S'}
\cM_h\qquad\forall n,m\in\ZZ.\label{egam}
\end{equation}
The elements of the subspaces $\cM_n$ are called
\textit{homogeneous elements of degree~$n$}.
\end{pf}

In the following, we sometimes refer to the constants $R$,
$S$, $R'$, and~$S'$ as the coefficients of diversity of the
grading.

For Krichever--Novikov algebras, the space of homogeneous
elements of degree $n$ for $n>0$ and $n<-g$ is
$\mathfrak{g}\otimes\cA_n$, where $\cA_n\subset\cA$ is the
subspace of functions that are of order $n$ at the point
$P_+$ and of order $-n+g$ at the point $P_-$ (cf. also
Example 2.1 below). For the constants $R$ and $S$ and other
details, see \cite{rKNFa,rSSkz}.

\subhead{(b) Holomorphic bundles. Tyurin parameters} Let
$\Sigma$ be a Riemann surface of genus $g$ and $F$ a
holomorphic bundle of rank $l$ over $\Sigma$. Suppose that
$l$ holomorphic sections $\Psi_1,\dots,\Psi_l$ of~$F$ are
given that form a basis in the fiber over any point except
for finitely many points $\ga_1,\dots,\ga_{gl}$, referred to
as \textit{points of degeneration}. This set of sections is
called a \textit{framing}, and a bundle that admits a
framing is said to be \textit{framed} \cite{rKNFd,rKNU}.

If one chooses a local trivialization of $F$ and defines
sections $\Psi_1,\dots,\Psi_l$ in this trivialization, say,
by column vectors (of functions), then one can take a matrix
$\Psi$ formed by these column vectors, which is
nondegenerate everywhere except for the points $\ga_i$,
$i=1,\dots,gl$, at which $\det\Psi$ has simple zeros for any
generic bundle, i.e., $\det\Psi(\ga_i)=0$ and
$(\det\Psi)'(\ga_i)\ne0$. In what follows, when speaking of
generic bundles, we mean just the validity of these
conditions. The divisor $D=\ga_1+\dots+\ga_{gl}$ is referred
to as the \textit{divisor of the bundle} $F$. Under our
conditions one has $\rank \Psi(\ga_i)=l-1$. Indeed,
$(\det\Psi)'(\ga_i)$ is the sum of determinants in which
$l-1$ rows are taken from the matrix $\Psi$ and one row from
the matrix $\Psi'$. By assumption, at least one of these
determinants is nonzero, and hence $\Psi$ has $l-1$ linearly
independent rows.

In this case, for any $i=1,\dots,gl$ there exists a
nontrivial solution of the system of linear equations
$\Psi(\ga_i)\a_i=0$; this solution is unique up to
proportionality. Let $\a_i=(\a_{ij})^{j=1,\dots,l}$
($i=1,\dots,gl$). In~\cite{rKNFd,rKNU}, the divisor of the
bundle and the numbers
$(\a_{ij})^{j=1,\dots,l}_{i=1,\dots,gl}$ are called the
\textit{Tyurin parameters} of the bundle~$F$. A framing is
defined uniquely up to a \textit{right} action of elements
of the group $GL(l)$ on the matrix $\Psi$ (while the
transition functions act on this matrix \textit{on the
left}; we see that the action of $GL(l)$ commutes with the
transition functions and hence takes sections to sections).
Owing to this fact, the Tyurin parameters are defined
uniquely up to proportionality and to multiplication by
elements of $GL(l)$ \textit{on the left}. According to the
results in \cite{rTur}, the Tyurin parameters define the
corresponding bundle uniquely up to equivalence. We point
out that the set of points $\ga_1,\dots,\ga_{gl}$ is the
same for equivalent framed bundles.

As was noted in \cite{rKNFd,rKNU}, Tyurin's results in
\cite{rTur} imply the following description of the space of
meromorphic sections of $F$ (we consider only meromorphic
sections that are holomorphic outside~$P_\pm$). The elements
$\Psi_j(P)$ ($j=1,\dots,l$) form a basis in the fiber over
any point $P$ outside the divisor $D$. Hence the value
$S(P)$ for any meromorphic section $S$ can be decomposed
with respect to this basis. To the section $S$ we assign the
vector function $f=(f_1,\dots,f_l)^T$ on the Riemann surface
$\Sigma$ in such a way that outside $D$ one has
\begin{equation}
S(P)=\sum_{j=1}^l\Psi_j(P)f_j(P).\label{razl}
\end{equation}
It follows from the Kronecker formula that $f_j=
\det(\Psi_1,\dots,\Psi_{j-1},S,\Psi_{j+1},\dots,\Psi_l)
(\det\Psi)^{-1}$, and therefore, the functions $f_j$ can be
extended to the points of the divisor $D$. The extensions
have at most simple poles at these points, since
$\Psi_1,\dots,\Psi_l$ are holomorphic at the points of $D$
and $\det\Psi$ has simple zeros at these points. In local
coordinates in a neighborhood of the point $\ga_i$, it
follows from \eqref{razl} that
$S(z)=\Psi(\ga_i)(\op{res}_{\ga_i}f)z^{-1}+ O(1)$. The
left-hand side of the last relation is holomorphic. Hence,
the residues of the functions $f_j$, $j=1,\dots,l$, at the
point $\ga_i$ satisfy the system of linear equations
$\Psi(\ga_i)(\op{res}_{\ga_i}f)=0$, which is just the system
satisfied by the Tyurin parameters at this point. By
assumption, the matrix $\Psi(\ga_i)$ is of rank $l-1$.
Hence, the vectors $\op{res}_{\ga_i}f$ and $\a_i$ are
proportional.

The following assertion represents the results of
\cite{rTur}.

\begin{tm}{Proposition 2.1 \cite{rKNFd,rKNU}}
For a generic bundle $F$, the space of meromorphic sections
holomorphic outside the distinguished points $P_\pm$ is
isomorphic to the space of meromorphic vector functions
$f=(f_1,\dots,f_l)^T$ on the same Riemann surface such that
$f$ is holomorphic outside $D$ and $P_\pm$, has at most
simple poles at the points of~$D$, and satisfies the
relations
$$
(\op{res}_{\ga_i}f_j)\a_{ik}
=(\op{res}_{\ga_{i}}f_k)\a_{ij},\qquad
i=1,\dots,gl,\;j=1,\dots,l.
$$
\end{tm}

\subhead{(c) Krichever--Novikov bases} In this subsection,
the exposition up to and including Proposition~2.3 follows
the note \cite{rKNRm}.

For any pair of integers $n$, $j$ ($0\le j<l$), we shall
construct a vector function in the space defined in
Proposition~2.1 and denote this function by~$\psi_{n,j}$.
The number $n$ is called the \textit{degree}
of~$\psi_{n,j}$. The function $\psi_{n,j}$ will be specified
by its asymptotic behavior at the points $P_{\pm}$. Let us
treat $\psi_{n,j}$ as a column vector and form a matrix
$\Psi_n$ of these vectors. At the point $P_+$, we set
\begin{equation}
\Psi_n(z_+)=z_+^n\sum_{s=0}^\infty\xi^+_{n,s}z_+^s,\quad
\text{where }\,\xi^+_{n,0}=\begin{pmatrix}
1&*&\hdots&*\\
0&1&\hdots&*\\
\vdots&\vdots&\ddots&\vdots\\
0&0&\hdots&1
\end{pmatrix},\label{bas1}
\end{equation}
and at the point $P_-$, we set
\begin{equation}
\Psi_n(z_-)=z_-^{-n}\sum_{s=0}^\infty\xi^-_{n,s}z_-^s,\quad
\text{where }\,\xi^-_{n,0}=\begin{pmatrix}
*&0&\hdots&0\\
*&*&\hdots&0\\
\vdots&\vdots&\ddots&\vdots\\
*&*&\hdots&*
\end{pmatrix},\label{bas2}
\end{equation}
where $z_\pm$ are the local parameters in neighborhoods of
the points $P_\pm$, respectively, and the asterisks $*$
stand for arbitrary complex numbers.

Thus, the function $\Psi_n$ has a zero of order $n$ at the
point $P_+$ and a pole of order $n$ at the point $P_-$, and
the determinant of $\Psi_n$ has $gl$ simple poles at the
points of the divisor $D$ and admits an additional divisor
of zeros outside the points $P_\pm$; the latter divisor is
not given \textit{a priori}.

\begin{pf}{Example 2.1}
In the Krichever--Novikov basis in the algebra
$\cA\,(\Sigma,P_\pm)$ with singularities at two points, the
asymptotics are of the following form \cite{rKNFa} (for $m>0$ and for
$m<-g$):
\begin{equation}
A_m\cong\a_m^{\pm}z_{\pm}^{\pm
m+\e_\pm}(1+O(z_\pm)),\qquad \a_m^+=1,\label{class}
\end{equation}
where $\e_+=0$ and $\e_-=-g$, i.e., the sum of orders (the
degree of the divisor of distinguished points) is equal to
$-g$ (note that $l=1$ in this case). Hence, there are
exactly $g$ zeros outside the points $P_\pm$. This example
differs from the general case studied above in that the
additional $g$ poles are concentrated at the same points
$P_\pm$ and the external divisor of degree $-g$ is absent.
This difference is unessential. The isomorphism is
established by the multiplication (division) by a scalar
function that has a zero of multiplicity $g$ at the point
$P_- $ and a given divisor of (simple) poles of degree $-g$
(such a function exists by the Riemann--Roch theorem).
\end{pf}

\begin{tm}{Proposition 2.2 \cite{rKNRm}}\quad
\rom1. There exists a unique matrix function $\Psi_n$
satisfying conditions \eqref{bas1} and \eqref{bas2}.

\rom2. The dimension of the spaces generated by the vector
functions $\psi_{n,j}$ for a given $n$ is equal to~$l$.
\end{tm}

\begin{pf}{Proof}
It is clear that assertions 1 and 2 are equivalent. In the
part related to the matrices $\xi^\pm_{n,0}$, conditions
\eqref{bas1} and \eqref{bas2} are the normalization
conditions that uniquely select the vector functions
$\psi_{n,j}$ in the $l$-dimensional space.

The dimension of the space of functions with $gl$ poles at
the points of the divisor $D$ with the orders $\pm n$ at the
points $P_\pm$, respectively, is equal to $gl-g+1$. Since
our vector functions are of dimension $l$, we find that the
dimension of the corresponding space is $l(gl-g+1)$.
However, these functions satisfy the $(l-1)gl$ Tyurin
relations. Therefore, the dimension of the function space
under consideration is equal to $l(gl-g+1)-(l-1)gl=l$. \qed
\end{pf}

The algebra $\cA\,(\Sigma,P_\pm)$ naturally acts both on the
sections of the bundle $F$ by the multiplication of these
sections by functions and on the corresponding vector
functions.

\begin{tm}{Proposition 2.3}
The action of the elements of the algebra
$\cA\,(\Sigma,P_\pm)$ on the basis elements $\psi_{n,j}$ is
quasigraded\rom:
\begin{equation}
A_m\psi_{n,j}
=\sum_{k=m+n}^{m+n+\bar{g}}\,\sum_{j'=0}^{l-1}
C^{k,j'}_{m,n,j}\psi_{k,j'},\label{act}
\end{equation}
where $\bar{g}=g+1$ for $-g\le m\le 0$ and $\bar{g}=g$
otherwise.
\end{tm}

\begin{pf}{Proof}
Let us use the asymptotic behavior \eqref{class}, where
$\e_-=-g$ for $m>0$ and for $m<-g$ and $\e_-=-g-1$ for
$-g\le m\le 0$ \cite{rKNFa}. Set $\bar{g}=-\e_-$.
Multiplying the matrix $\Psi_n$ by $A_m$, we obtain the
order $z_+^{n+m}$ at the point $P_+$ and the order
$z_-^{-n-m-\bar{g}}$ at the point $P_-$. This is sufficient
for the validity of relations \eqref{act} for some values of
the constants $C^{k,j'}_{m,n,j}$. In particular, the
diversity of the grading with respect to the index $n$ in
the general case is equal to $g$ (for $-g\le m\le 0$ this
index is equal to $g+1$), and hence it is bounded. This
proves that the action is quasigraded. In what follows we
use the fact that the diversity of the grading with respect
to the index $ln-j$ is equal to $lg$ in the general case.
\qed
\end{pf}

In contrast with the algebra $\cA\,(\Sigma,P_\pm)$, the Lie
algebra $\mathfrak{g}$ admits no natural action on sections
of the bundle. However, suppose that a representation $\tau$
of the Lie algebra $\mathfrak{g}$ in an $l$-dimensional
space is given. Then this representation induces an action
of this algebra on the vector functions~$\psi_{n,j}$.
Namely, let $x\in\mathfrak{g}$ and let $\tau(x)=(x^i_{i'})$,
where $(x^i_{i'})$ is a matrix and the indices $i$ and $i'$
range from $1$ to $l$. We represent the vector function
$\psi_{n,j}$ in the form of the set of its coordinates,
$\psi_{n,j}=(\psi_{n,j}^{i'})$. In the coordinate
representation, the action of $x$ on $\psi_{n,j}$ can be
written in the form
$(x\psi_{n,j})^i=\sum_{i'}x^i_{i'}\psi_{n,j}^{i'}$. Having
the action of the algebras $\mathfrak{g}$ and
$\cA\,(\Sigma,P_\pm)$ on vector functions, we define the
corresponding action of the Lie algebra
$\mathfrak{g}\otimes\cA\,(\Sigma,P_\pm)$ as follows:
\begin{equation}
((xA_m)\psi_{n,j})^i
=\sum_{k=m+n}^{m+n+\bar{g}}\sum_{j'=0}^{l-1}\sum_{i'}
C_{m,n,j}^{k,j'}x^i_{i'}\psi_{k,j'}^{i'}.
\label{actg}
\end{equation}

Formula \eqref{actg} is an expression for the tensor product
of the representation $\tau$ of $\mathfrak{g}$ and the
action \eq{2.9} of $\cA\,(\Sigma,P_\pm)$. Thus, this formula
defines a representation of the tensor product of these
algebras. It follows from Proposition~2.3 that the action
\eqref{actg} is quasigraded and the diversity of the grading
with respect to the index $n$ (respectively, $ln-j$) is
equal to $g$ (respectively, $gl$).

In the sequel, we interpret the $\psi^i_{n,j}$ as formal
symbols rather than coordinates of vector functions and
treat relations \eqref{actg} as formal linear
transformations of these symbols. Under this approach, the
action of $\mathfrak{g}$ does not change the indices $n$ and
$j$, and the role of a highest vector is played by the
symbol $\psi^l_{n,j}$ for given $n$ and $j$.

\section{Fermion Representations}

Let us consider a framed holomorphic bundle $F$ of rank $l$
and the corresponding space $F_{KN}$ of Krichever--Novikov
vector functions. In this space, we consider the basis
$\{\psi_{n,j}\mid n\in\ZZ,\,j=0,\dots,l-1\}$, the
corresponding set of symbols $\psi_{n,j}^i$ ($i=1,\dots,l$),
and the action \eqref{actg} on this set; these objects were
defined in the preceding section. Recall that this action
depends on some representation
$\tau\colon\mathfrak{g}\mapsto\mathfrak{gl}(l)$. Starting
from these data, we will construct a space $V_F$ of
semi-infinite forms and a representation $\pi_{F,\tau}$ of
$\hat{\mathfrak{g}}$ in this space.

Let us number the symbols $\psi_{n,j}^i$ in
lexicographically ascending order of the triples $(n,-j,i)$,
setting the index of the element $\psi_{0,0}^l$ to be equal
to $-1$. Let $N=N(n,j,i)$ be the index of a triple
$(n,j,i)$; then $N(n,j,i)=l^2n-lj+i-l-1$. Set
$\psi_N=\psi_{n,j}^i$. Formal finite linear combinations of
the symbols $\{\psi_N\mid N\in\ZZ\}$ form a vector space,
which we denote by $F_{KN}^s$.

The space $V_F$ is constructed in the following way. It is
spanned over $\CC$ by formal expressions
$\psi_{N_0}\wedge\psi_{N_1}\wedge\dots$, referred to as
semi-infinite monomials, where $N_0<N_1<\dots$\,. It is
required that, under a transposition of $\psi_N$ and
$\psi_{N'}$ ($N\ne N'$), a semi-infinite monomial changes
its sign and one has $N_k=k+m$ for sufficiently large $k$
($k\gg 1$), where $m$ is an integer; following \cite{rKaRa},
$m$ is called the \textit{charge} of the monomial (the
reader should not confuse this with the \textit{central
charge}: it follows from Remark 3.1 that the charge $m$ is
rather a component of the weight).

By the \textit{degree of a monomial} $\psi$ of charge $m$
one means
\begin{equation}
\deg\psi =\sum_{k=0}^\infty (N_k-k-m).
\label{deg}
\end{equation}
The degree thus defined equips the space $V_F$ with a
quasigraded structure, i.e., a decomposition into a direct
sum of finite-dimensional subspaces. Note that the
definition~of the numbers $N_k$ is ambiguous; namely, it
depends on the numbering of the elements $\psi_{n,j}^i$ for
a given $n$. It is of importance that the definition~of
degree is independent of this ambiguity since the
modification of the above numbering amounts to a permutation
of some summands in formula \eqref{deg}.

The representation $\pi_{F,\tau}$ is defined as follows.
Consider the Lie algebra  $\bar{\mathfrak{g}}=\mathfrak{g}
\otimes\cA\,(\Sigma,P_\pm)$ embedded in the Lie algebra
$\hat{\mathfrak{g}}$ as a linear space. Every element
$x\in\bar{\mathfrak{g}}$ acts on the symbols $\psi_N$ by
linear changes in a quasigraded way by formula \eqref{actg}.
This precisely means that if the symbols $\psi_N$ are
regarded as a formal basis of an infinite-dimensional linear
space, then the action of an element
$x\in\bar{\mathfrak{g}}$ in this basis is defined by an
infinite matrix with finitely many nonzero diagonals, or,
equivalently, by a difference operator. Following
\cite{rKac}, we denote the algebra of these matrices
by~$\mathfrak{a}_{\infty}$. Thus, to the basis $\{\psi_N\}$
there corresponds an embedding of $\bar{\mathfrak{g}}$ in
the algebra $\mathfrak{a}_{\infty}$. It suffices to define
the action of the algebra $\mathfrak{a}_{\infty}$ on $V_F$,
and then, by virtue of this embedding, we obtain a
representation of $\hat{\mathfrak{g}}$ in the same space.

Let us take a basis element
$E_{IJ}\in\mathfrak{a}_{\infty}$. We define its action on a
semi-infinite monomial
$\psi=\psi_{I_0}\wedge\psi_{I_1}\wedge\dots$ as follows:
\begin{equation}
r(E_{IJ})\psi=(E_{IJ}\psi_{I_0})\wedge\psi_{I_1}\wedge\ldots+
\psi_{I_0}\wedge(E_{IJ}\psi_{I_1})\wedge\ldots+\dots\,.
\label{infact}
\end{equation}
By virtue of the condition $I_k=k+m$ ($k\gg 1$), the action
\eqref{infact} is well defined for $I>J$ and for $I<\nobreak
J$. The operator $r(E_{II})$ acts on a semi-infinite
monomial by multiplication by $1$ or by $0$, depending on
whether the monomial contains $\psi_I$. Therefore, for a
diagonal matrix $D=\sum_{i=-\infty}^\infty \la_IE_{II}$ with
infinitely many coefficients $\la_I\ne 0$ such that $I>0$,
the sum $D\psi$ (which is well defined by \eqref{infact})
can contain infinitely many terms. In this case, we apply
the following regularization \cite{rKaRa}. Let
\begin{alignat}2
\hat{r}(E_{IJ})&=r(E_{IJ})&&\quad\text{for }\,I\ne J\,\text{
and for }\,I=J<0,
\label{reg1}\\
\hat{r}(E_{II})&=r(E_{II})-\op{Id}&&\quad\text{for }\,I\ge 0.
\label{reg2}
\end{alignat}
Then the following commutation relations hold \cite{rKaRa}:
\begin{align}
[\hat{r}(E_{IJ}),\hat{r}(E_{MN})]&=\delta_{{}_{JM}}{\hat
r}(E_{IN})-\delta_{{}_{NI}}\hat{r}(E_{MJ})\quad\text{for
}\,(I,J)\ne(N,M),\label{com1}\\
[\hat{r}(E_{IJ}),\hat{r}(E_{JI})]&
=\hat{r}(E_{II})-\hat{r}(E_{JJ})+\op{Id}.
\label{com2}
\end{align}
Thus, the regularization gives the following cocycle $\a$:
\begin{alignat*}2
\a(E_{IJ},E_{JI})&=-\a(E_{JI},E_{IJ})=1&&\quad\text{for
}\,I\ge0,\;J<0,\\
\a(E_{IJ},E_{MN})&=0&&\quad\text{otherwise}.
\end{alignat*}

It is quite clear that the action of the algebra
$\mathfrak{a}_{\infty}$ preserves the charge of a monomial,
since the infinite ``tail'' of each of the monomials on the
right-hand side in \eqref{infact} coincides with the
``tail'' of the monomial on the left-hand side in
\eqref{infact}. Moreover, $\pi_{F,\tau}$ is a quasigraded
representation. This follows from the fact that the action
on the space $F_{KN}$ is quasigraded (Proposition 2.3 and
relations \eqref{act} and \eqref{actg}).

By \textit{highest} (or \textit{vacuum}) \textit{monomials}
we mean semi-infinite monomials of the form $\tilde\psi_M
=\psi_M\wedge\psi_{M+1}\wedge\psi_{M+2}\wedge\dots$ (the
indices are successive, starting from some number).

The algebra $\cA\,(\Sigma,P_\pm)$ can also be embedded in
the algebra $\mathfrak{a}_{\infty}$ of matrices with
finitely many nonzero diagonals. To prove this fact, it
suffices to rewrite relation \eqref{act} by replacing every
$\psi_{n,j}$ by $\psi_{n,j}^i$ and every triple of indices
$n$, $j$, $i$ by the corresponding index $N=N(n,j,i)$:
$$
A_m\psi_N=\sum_{N'=l^2m+N}^{l^2m+N+g_l}
C^{N'}_{m,N}\psi_{N'},
$$
where $g_l=l^2\bar{g}$ (see Proposition 2.3). The following
lemma clarifies the character of the action of the algebra
$\cA\,(\Sigma, P_\pm)$ on the highest monomials. Let us take
a highest monomial
$\tilde\psi_M=\bigwedge_{N=M}^\infty\psi_N$.

\begin{tm}{Lemma~3.1}
The element $A_m\tilde\psi_M$ has a nonzero projection on
$\tilde\psi_M$ only for $m=-g,\dots,-1,0$, and this
projection is equal to $\sum_{K=M}^{-1}C_{mK}^K$.
\end{tm}

\begin{pf}{Proof}
It follows from relation \eqref{act} that, for arbitrary
integer $N$, the element $A_j\psi_N$ has a nonzero
projection on $\psi_N$ if and only if $\psi_N$ occurs on the
right-hand side in \eqref{act}, i.e., only for $n+m\le n\le
n+m+g$. This implies that $-g\le m\le 0$.

Let us represent the action of an element $A_m$, $-g\le
m\le0$, by an infinite matrix with finitely many nonzero
diagonals. Then one has $A_m=\dots+\sum_K
C_{mK}^KE_{KK}+\dots$, where the dots stand for terms with
$E_{IJ}$ for $I\ne J$. For $I<J$, these terms send
$\tilde\psi_M$ to $0$ and, for $I>J$, to either $0$ or a
monomial of some degree less than that of $\tilde\psi_M$.
For $E_{KK}$, where $K$ is arbitrary, one has
$E_{KK}\psi_N=\delta_{{}_{KN}}\psi_N$, and for $K\ge 0$ it
follows from the above regularization that
$E_{KK}\tilde\psi_M=0$. Hence,
\begin{equation}
A_m\tilde\psi_M=\bigg(\sum_{K=M}^{-1}C_{mK}^K\bigg)
\tilde\psi_M+\dots,\label{wght}
\end{equation}
where the dots stand for the sum of all terms whose degrees
are less than that of $\tilde\psi_M$.\qed
\end{pf}

\begin{pf}{Remark~3.1}
It follows from the proof of Lemma~3.1 that the eigenvalue
of the operator $\op{Id}=\sum_NE_{NN}$ (the identity
operator of the space $F_{KN}^s$) on $\tilde\psi_M$
coincides with $-M$ for $M<0$ and is zero for $M\ge0$.
\end{pf}

\begin{pf}{Definition~3.1}
Two fermion representations are said to be \textit{strongly
equivalent \rom(isomorphic\rom)} if there exists an
isomorphism of the underlying linear spaces that commutes
with the action of the algebra $\hat{\mathfrak{g}}$, is
\textit{quasi-homogeneous} (i.e., takes every homogeneous
component of the first module to a homogeneous component of
the other), and takes any highest monomial to a highest
monomial of the same charge.
\end{pf}

\begin{tm}{Theorem~3.1}
Let a fermion representation $\pi$ be defined by a framed
holomorphic bundle $F$ and by an irreducible representation
$\tau$ of the algebra $\mathfrak{g}$, and let a fermion
representation $\pi'$ be defined by a framed holomorphic
bundle $F'$ and by an irreducible representation $\tau'$.
The representations $\pi$ and $\pi'$ are equivalent if and
only if the corresponding framed bundles $F$ and $F'$ are
equivalent and the representations $\tau$ and $\tau'$ of the
algebra $\mathfrak{g}$ are equal.
\end{tm}

\begin{pf}{Proof}
Let the bundles $F$ and $F'$ be equivalent, and let the
equivalence be defined by a mapping $C\colon F\to F'$.
Obviously, $C$ induces the following isomorphisms:

\eq{1} of the spaces of global holomorphic sections,
$H^0(\Sigma,F)\cong H^0(\Sigma,F')$; under this isomorphism,
the matrices $\Psi$ and $\Psi'$ composed of the basis
holomorphic sections correspond to each other,
$\Psi'=C\Psi$;

\eq{2} of the spaces of meromorphic Krichever--Novikov
vector functions regarded as $\cA\,(\Sigma,P_\pm)$-modules:
$F_{KN}\cong F'_{KN}$, and the same holds for the related
spaces $F_{KN}^s$ and $F_{KN}^{\prime s}$.

Since the representations $\tau$ and $\tau'$ are equal, it
follows that the spaces $F_{KN}^s$ and $F_{KN}^{\prime s}$
are isomorphic as $\mathfrak{g}$-modules as well. Hence, the
corresponding representations of the algebra
$\hat{\mathfrak{g}}$ in the spaces of semi-infinite forms
over $F_{KN}^s$ and $F_{KN}^{\prime s}$ are equivalent.

Let us pass to the proof of the second part of the theorem.
Consider two equivalent fermion representations $\pi$ and
$\pi' $ of $\hat{\mathfrak{g}}$ generated by framed bundles
$F$ and $F'$ and by irreducible representations $\tau$ and
$\tau'$, respectively, in the sense of the above
construction. Let $C$ be a strong isomorphism, $\pi'=C\pi
C^{-1}$.

Let $v$ be a highest monomial in the representation space of
$\pi$, and let $v'=Cv$. Then it follows from Definition 3.1
that $v'$ is also a highest monomial. Our immediate aim is
to prove that the representations $\tau$ and $\tau'$ are
equivalent. To this end, we consider the structure of the
highest monomials $v$ and $v'$.

Let $v=\tilde\psi_M$. A triple of indices $(n,j,i)$ such
that $M=N(n,j,i)$ is determined uniquely (see above). We
write $n=n(M)$, $j=j(M)$, and $i=i(M)$. Let us choose $v$ in
such a way that $i=l$. As was noted at the end of the
previous section, this is equivalent to the condition that
$\psi_M$ is a highest vector of a subrepresentation
equivalent to $\tau$ and acting on the space $F_{KN}^s$. Let
$v'=\tilde\psi'_{M'}$, and accordingly, let $n'=n(M')$,
$j'=j(M')$, and $i'=i(M')$. By definition, a
quasi-homogeneous isomorphism preserves the charge.
Therefore, $M=M'$, and hence $i=i'=l$ and $\psi'_{M'}$ is a
highest vector of a subrepresentation equivalent to $\tau'$.

Thus, $\psi_M$ and $\psi'_{M'}$ are highest vectors of the
representations $\tau$ and $\tau'$, respectively. Hence, the
structure of the monomial $v=\tilde\psi_{M}$ is as follows:
the first place is occupied by the highest vector
${\psi}_{M}$ of the representation $\tau$, and to the right
of this place one can find an infinite exterior product of
the symbols $\psi_N$, which can be partitioned into
consecutive groups (we call them bags) each of which is the
exterior product of (all) basis elements of the
representation $\tau$. This readily implies that
\begin{equation}
\pi
(x)v=(\tau(x)\psi_M)\wedge\psi_{M+1}\wedge
\dots\label{deistv1}
\end{equation}
 for an arbitrary $x\in\mathfrak{g}$. Indeed, if $\tau(x)$
is a nilpotent element, then the other terms (that occur
under the action of $\pi(x)$ on $v$ by the Leibniz formula)
are equal to $0$. If $\tau(x)$ is a diagonal element, then
their sum is annihilated for any bag, since the sum of the
weights (counted according to their multiplicity) of a
finite-dimensional representation of a semisimple Lie
algebra is equal to $0$. (This is obvious for any
irreducible representation of $sl(2)$, and the other cases
can be reduced to this one by decomposing into irreducible
representations of simple three-dimensional subalgebras
corresponding to simple roots. In the case of a reductive
algebra, this argument would fail for a central element.) A
similar consideration shows that
\begin{equation}
\pi' (x)v'=(\tau'(x)\psi'_{M'})\wedge\psi'_{M'+1}\wedge\dots\,.
\label{deistv2}
\end{equation}
Let $x=h$ be a Cartan element, and let $\a_\tau$ and
$\a_{\tau'}$ be the highest weights of the representations
$\tau$ and $\tau'$, respectively. Then it follows from
formulas \eqref{deistv1} and \eqref{deistv2} that
$\pi(h)v=\a_\tau(h)v$ and $\pi'(h)v'=\a_{\tau'}(h)v'$. Since
the representations $\pi$ and $\pi'$ are equivalent, it
follows that $\a_\tau=\a_{\tau'}$. Hence, $\tau$ and $\tau'$
are equivalent.

The Lie algebra $\mathfrak{g}$ acts on the spaces $F_{KN}^s$
and $V_F$, and hence the same holds for the universal
enveloping algebra $U(\mathfrak{g})$. The associative
algebra $\cA\,(\Sigma, P_\pm)$ and the tensor product
$U(\mathfrak{g})\otimes\cA\,(\Sigma, P_\pm)$ also act on
these spaces. Moreover, by construction, the representation
of the algebra $U(\mathfrak{g})\otimes\cA\,(\Sigma, P_\pm)$
is decomposed into a product of representations of the
algebras $U(\mathfrak{g})$ and $\cA\,(\Sigma,P_\pm)$.
Namely, if $u\otimes A\in
U(\mathfrak{g})\otimes\cA\,(\Sigma,P_\pm)$, then
$\pi(u\otimes A)=\tau(u)\pi(A)$. We claim that the
equivalence of representations of the algebra
$\hat{\mathfrak{g}}$ implies the equivalence of the
corresponding representations of the algebra
$\cA\,(\Sigma,P_\pm)$. Let $u=u^{i_1\dots i_n}e_{i_1}\dots
e_{i_n}\in U(\mathfrak{g})$ (where summation with respect to
repeated sub- and superscript is assumed). Consider the
element $u'=u^{i_1\dots i_n}e_{i_1}\dots(e_{i_n}A)$ in which
the function $A$ in any monomial is placed in the last
factor. Then $\pi(u')=\tau(u)\pi(A)$. For the equivalent
representation $\pi'$ of the algebra $\hat{\mathfrak{g}}$,
one has $\pi' =C\pi C^{-1}$, where $C$ is an intertwining
operator, and hence
\begin{equation}
\tau'(u)\pi' (A)=C(\tau(u)\pi (A))C^{-1}. \label{equiv1}
\end{equation}
Let $Z(\mathfrak{g})$ be the center of the ring
$U(\mathfrak{g})$. Since the representations $\tau$ and
$\tau'$ are irreducible, it follows that the ring
$Z(\mathfrak{g})$ acts on $F_{KN}^s$ and $F_{KN}^{\prime s}$
by scalar operators. For $u\in Z(\mathfrak{g})$, we set
$\tau(u)=\la_u\circ\op{Id}$ and
$\tau'(u)=\la_u'\circ\op{Id}$. Since $\tau\cong\tau'$, it
follows that $\la_u=\la_u'$. Now it follows from relation
\eqref{equiv1} that, cancelling $\la_u$ (which can be done
for at least one element $u$), we obtain
$\pi'(A)=C\pi(A)C^{-1}$, as desired.

Let us show that the representations of the ring
$\cA\,(\Sigma,P_\pm)$ in the spaces $F_{KN}$ and $F'_{KN}$
are equivalent. To this end, consider the action of an
element $A_n\in\cA\,(\Sigma,P_\pm)$ on monomials of the form
$\psi_{N,M}=\psi_N\wedge\psi_{M+1}\psi_{M+2}\dots$. If $M-
N$ is large compared with $n$, then the structure constants
of the action of the element $A_n$ on $\psi_N$ are contained
in the set of the structure constants of the action of $A_n$
on $\psi_{N,M}$. Thus, one can recover the first action from
the other one. As usual, let $N=N(n,j,i)$. The structure
constants of the action of the element $A_m$ on
$\psi_{n,j}^i$ are independent of $i$ and coincide with the
structure constants of the action of $A_m$ on $\psi_{n,j}$,
which can  be observed by comparing formulas \eqref{act} and
\eqref{actg}. Thus, we have obtained the desired equivalence
of $F_{KN}$ and $F'_{KN}$ treated as
$\cA\,(\Sigma,P_\pm)$-modules. As was shown in \cite{rKNRm},
in this case the framed bundles $F$ and $F'$ are equivalent.
Let us return to the representations $\tau$ and $\tau'$. As
was shown above, they are equivalent, i.e., there exists an
element $\ga\in GL(l)$ such that $\tau'=\ga\tau\ga^{-1}$. If
the matrix $\ga$ were not scalar, then it would ``mix''
homogeneous components $\psi^{\prime i}_{n,j}$ with respect
to the index $i$. Thus, the isomorphism $C$ would be not
quasi-homogeneous, which contradicts our assumption. Hence,
$\ga$ belongs to the center of the group $GL(l)$, and
therefore, $\tau=\tau'$. This completes the proof. \qed
\end{pf}

There is 
a more general notion of equivalence for fermion
representations than that in Definition~3.1.

\begin{pf}{Definition~3.2}
Two fermion representations are said to be
\textit{equivalent \rom(isomorphic\rom)} if there is an
isomorphism of the form $C_\ga=C\cdot{\tilde\ga}$ between
the representation spaces that commutes with the action of
the algebra $\hat{\mathfrak{g}}$, where $C$ is a strong
isomorphism in the sense of Definition~3.1, $\ga\in GL(l)$,
and the mapping $\tilde\ga$ is defined as follows:
${\tilde\ga}(\psi_{N_1}\wedge\psi_{N_2}\wedge\dots)=
\ga\psi_{N_1}\wedge\ga\psi_{N_2}\wedge\dots$.
\end{pf}

Each of the mappings $C$ and $\tilde\ga$ in this definition
appears as a result of some operation over bundles. As was
shown in the proof of Theorem~3.1, to a strong isomorphism
$C$ there corresponds a bundle isomorphism that can be
defined by an equivalent change of the gluing functions. In
Theorem~3.2 we will prove that the mapping $\tilde\ga$
arises from the change of framing in a bundle. Any two
operations of these two types commute (see \S2,
Sec.~(b)). Hence, the mappings $C$ and $\tilde\ga$
commute, and thus an isomorphism of the form $C_\ga$ is an
equivalence relation indeed. The following lemma~will be
used in the proof of Theorem~3.2.

\begin{tm}{Lemma~3.2}
For any irreducible $l$-dimensional representation $\tau$ of
the Lie algebra $\mathfrak{g}$ and any generic holomorphic
bundle $F$ of rank $l$ and of degree $gl$, the equivalence
class of the representation $\tau$ is in a one-to-one
correspondence with the set of framings of the bundle $F$.
\end{tm}

\begin{pf}{Proof}
The group $GL(l)$ acts on the irreducible representations of
the algebra $\mathfrak{g}$ in $\CC^l$ as follows: an element
$\ga\in GL(l)$ takes the representation $\tau$ to the
equivalent representation $\ga^{-1}\tau\ga$, and then the
representation $\tau$ remains stable if and only if the
element $\ga$ is scalar, $\ga\in\CC\cdot\op{Id}$.

At the same time, the group $GL(l)$ acts on the framings of
the bundle $F$ as follows: an element $\ga$ takes a framing
defined by the matrix $\Psi$ of sections to the framing
defined by the matrix $\Psi\ga$ of sections (see \S2,
Sec.~(b)). Moreover, the scalar operators (and only these
operators) preserve the Tyurin parameters as points of the
projective space, since they multiply all parameters by the
same number.

Thus, both sets under consideration are equal to
$GL(l)/\CC\cdot\op{Id}$ as manifolds.\qed
\end{pf}

\begin{tm}{Theorem~3.2}
The isomorphism relation introduced by Definition~\rom{3.2}
is an equivalence relation on the set of fermion
representations. The corresponding set of equivalence
classes of fermion representations is in a one-to-one
correspondence with the set of pairs $([F],[\tau])$, where
$F$ is an equivalence class of \rom(nonframed\rom)
holomorphic bundles of rank $l$ and of degree $gl$ \rom(in
general position\rom) and $[\tau]$ is an equivalence class
of\/ $l$-dimensional irreducible representations of the
algebra~$\mathfrak{g}$.
\end{tm}

\begin{pf}{Proof}
Let a generic holomorphic (nonframed) bundle $F$ of rank $l$
of degree $gl$ and an irreducible $l$-dimensional
representation $\tau$ of the algebra $\mathfrak{g}$ be
given. Choose any two framings in $F$ that are defined by
matrices of holomorphic sections $\Psi$ and $\Psi'$. Then
there exists a matrix $\ga\in GL(l)$ such that
$\Psi'=\Psi\ga^{-1}$. Suppose that to the data
$\{F,\tau,\Psi\}$ ($\{F,\tau',\Psi'\}$) there corresponds a
fermion representation  $\pi=\pi(F,\tau,\Psi)$
($\pi'=\pi(F,\tau',\Psi')$, respectively). In this case, we
have an isomorphism $\tilde\ga$ between the representation
spaces $\pi$ and $\pi'$: by the construction of the fermion
representation, for each monomial
$\psi_{N_1}\wedge\psi_{N_2}\wedge\dots$ (if it belongs to
the representation space of~$\pi$) we have the corresponding
monomial $\ga\psi_{N_1}\wedge\ga\psi_{N_2}\wedge\dots$. This
can be seen from relation \eqref{razl}. For the framing
$\Psi$, this relation (in matrix form) becomes $S=\Psi f$.
Respectively, for the framing $\Psi'$ this relation is of
the form $S=\Psi\ga^{-1}\cdot\ga f$, i.e., the elements
$f\in F_{KN}$ are subjected to the automorphism $\ga$.

Thus, we have proved that the mapping $\tilde\ga$ appears as
the result of the change of framing
$\Psi\mapsto\Psi\ga^{-1}$ in the holomorphic bundle. As was
noted above, this implies that an isomorphism in the sense
of Definition~3.2 is an equivalence relation on the set of
fermion representations.

For the mapping $\tilde\ga$ to commute with the action of
the algebra $\mathfrak{g}$ (this is also one of the
conditions in Definition~3.2), it is necessary and
sufficient that $\tau'=\ga\tau\ga^{-1}$. Thus, we obtain the
following picture. By Theorem~3.1, the classes of strong
equivalence of fermion representations are parametrized by
the data sets $\{F,\tau,\Psi\}$. For a given structure of a
holomorphic bundle $F$, to an isomorphism $\tilde\ga$ there
corresponds a transformation of the related data set,
$\Psi\mapsto\Psi\ga^{-1}$ and $\tau\mapsto\ga\tau\ga^{-1}$.
Hence, to the \textit{equivalence} classes of fermion
representations (in the sense of Definition~3.2) we can
assign the orbits in the space of the data sets
$\{F,\tau,\Psi\}$ under the action of the transformations
$F\mapsto F$, $\Psi\mapsto\Psi\ga^{-1}$, and
$\tau\mapsto\ga\tau\ga^{-1}$. It follows from Lemma~3.2 that
these orbits are in a one-to-one correspondence with the
pairs $([F],[\tau])$.\qed
\end{pf}


\end{document}